\documentclass{article}
\usepackage[T1]{fontenc}
\usepackage{amsmath,amssymb}
\usepackage{graphicx,xcolor}
\usepackage{hyperref}
\newtheorem{remark}{Remark}
\title{CmcMesh: a python implementation of the DPW method}
\author{Thomas Raujouan, Martin Traizet}
\begin{document}
\maketitle
\begin{abstract}
	Minimal and constant mean curvature (CMC) surfaces in three-dimensional	space forms can be constructed with the Dorfmeister--Pedit--Wu (DPW) method. We introduce CmcMesh, a Python implementation covering the whole	pipeline from holomorphic potential to rendered surfaces, using Hoffman and Hoffman's Mesh algorithm. 
	We implement the Weierstrass representation for minimal surfaces in $\mathbb{R}^3$, Bryant's representation for CMC-1 surfaces in $\mathbb{H}^3$, and the DPW method for CMC surfaces in $\mathbb{R}^3$ and for minimal and CMC surfaces	in $\mathbb{S}^3$. The package is modular, making it easy to add new DPW variants without modifying existing code.
\end{abstract}
\section{Introduction}
The CmcMesh project is intended for differential geometers working in the field of constant mean curvature (CMC) surfaces, who are familiar with the DPW method ~\cite{dorfmeister1998weierstrass}.
It combines the Mesh algorithm ~\cite{callahan1988computer} with the DPW method.
Its main purpose is to produce pictures or make numerical experiments with the DPW method.
\medskip

The Mesh algorithm was developed in the 80's by David Hoffman and Jim Hoffman to make pictures of minimal surfaces in euclidean space $\mathbb R^3$ via the classical Weierstrass Representation.
It is a simple and robust algorithm which simultaneously triangulates the domain and computes the minimal immersion. The algorithm adapts the size of the triangles in the domain to the conformal factor of the immersion and its curvature, so that the triangulated surface looks nice.
\medskip

The DPW method is a loop-group representation for CMC surfaces in euclidean space $\mathbb R^3$ which is quite more elaborate and computationally demanding than Weierstrass Representation, but retains the idea that the surface can be constructed from some holomorphic data on a Riemann surface. The DPW method has many avatars so can also be used to compute CMC surfaces in the sphere $\mathbb S^3$, hyperbolic space
$\mathbb H^3$ and other spaces. More generally, it is a method to construct harmonic maps from a Riemann surface to a symmetric space.
\medskip

The DPW method has been implemented by Nick Schmitt in Cmclab and Xlab.
Over the years, Nick has produced hundreds of beautiful pictures of CMC surfaces and has been very helpful to the DPW community. Unfortunately for us, Nick retired a couple of years ago. The CmcMesh project started from the need to produce illustrations for some recent papers.
\medskip

CmcMesh is coded in python. All loop-group computations are fully vectorized in numpy for computational efficiency, including Iwasawa decomposition.

\medskip

The project includes a viewer to visualize the triangulated domain and surface, and methods to export to standard formats for 3D objects.
\section{How to obtain CmcMesh}
To get started, please download CmcMesh from the gitlab repository below and follow the instructions:
\begin{center}
\href{https://scm.univ-tours.fr/projetspublics/idp/cmcmesh}{scm.univ-tours.fr/projetspublics/idp/cmcmesh}
\end{center}
CmcMesh is intended to be a collaborative project and we would greatly appreciate that you send us any feedback and bug reports at \href{mailto:martin.traizet@univ-tours.fr}{martin.traizet@univ-tours.fr}
or \href{mailto:thomas.raujouan@univ-tours.fr}{thomas.raujouan@univ-tours.fr}.
Please let us know if you want to contribute the project by adding new examples or avatars of the DPW method.
\section{Architecture of the project}
As already said, the DPW method has several avatars. This is reflected in the architecture of the project.
\begin{center}
\includegraphics[width=8cm]{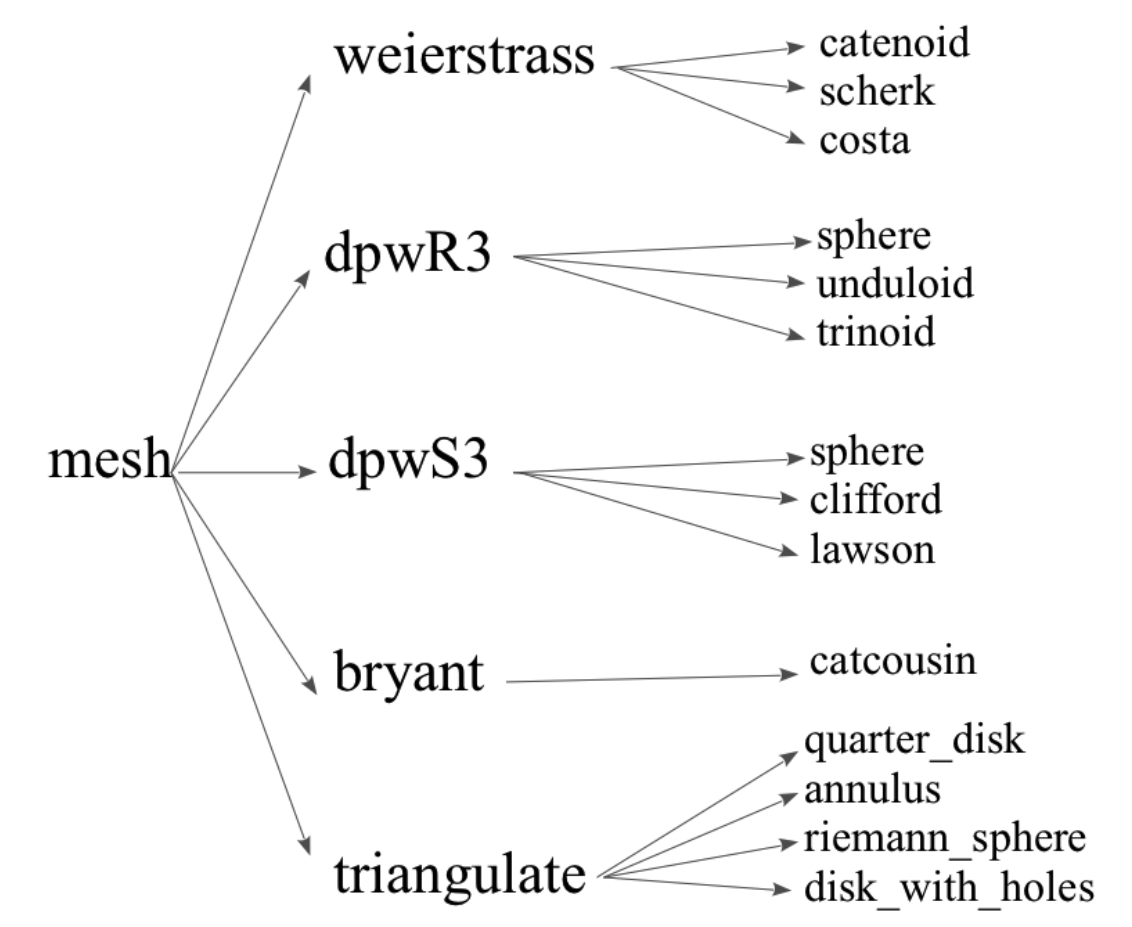}
\end{center}
The main module is \verb$mesh$ which implements the Mesh algorithm. The main class in the
module \verb$mesh$ is \verb$domain$, which allows the user to specify the domain on which the surface is parametrized (see Section \ref{section:domain}).

Then come several modules which implement the various avatars of the DPW method.
In each of these modules, the main class is \verb$data$, which allows the user to specify the holomorphic data for his surface (see Section \ref{section:data}). For example:
\begin{itemize}
\item \verb$dpwR3$ implements the DPW method for CMC-1 surfaces in $\mathbb R^3$. The holomorphic data here is the DPW potential $\xi$ and the initial condition $\Phi_0$.
\item \verb$dpwS3$ implements the DPW method for minimal and CMC surfaces in $\mathbb S^3$. The data is the DPW potential $\xi$, the initial condition $\Phi_0$ and the two Sym-points
$\lambda_1$, $\lambda_2$.
\item \verb$weierstrass$ implements the Weierstrass Representation for minimal surfaces in $\mathbb R^3$.
Here the data is the Gauss map $g$ and the holomorphic 1-form $\omega$.
\item \verb$bryant$ implements Bryant's Representation for CMC-1 surfaces in $\mathbb H^3$.
\item Finally, \verb$triangulate$ can be used to just triangulate a domain in the complex plane. Here the data consists of a single real-valued function $f(z)$ which tells the desired edge length, depending on $z$. This is here mainly to test and illustrate functionalities of the Mesh algorithm.
\end{itemize}
Then come examples.
A typical example file contains:
\begin{itemize}
\item the definition of the domain,
\item the definition of the data,
\item a call to the \verb$grow$ method which will start the Mesh algorithm and return a triangulated surface in the appropriate space,
\item various post-processing operations, such as: extending the surface by symmetry, extracting curvature lines, exporting...
\item Finally, a call to the method \verb$view$ to launch the viewer.
\end{itemize}
As a user of the CmcMesh project, if your favorite avatar of the DPW method is already implemented, all you have to do is to start from an existing example and adapt it to suit your needs.
We plan to implement more avatars of the DPW method, but if your favorite avatar is not yet implemented, you should be able to do it yourself. For this, you will need to understand how the avatar interacts with the Mesh algorithm: see Section \ref{section:how}.
\medskip

The main reason to have a clean separation between the Mesh algorithm and the various avatars is that new functionalities of the Mesh algorithm can be developed without having to change a single line of code in the avatars.

\section{Structure of an example file}
\label{section:structure}
\subsection{Defining the domain}
\label{section:domain}
The domain is defined by the following lines:
\begin{verbatim}
from lib.mesh import domain
mydomain=domain(
    equations=[f1,f2,...],
    cuts=[c1,c2,...],
    origin=z0,
    maxNbPoints=1000
)
\end{verbatim}
\begin{itemize}
\item Each $f_i$ is a real-valued function of the complex variable. The domain is defined by
$f_i(z)>0$ for all $i$.
For example, $f(z)=R-|z-a|$ defines the disk $D(a,R)$. Alternately, one can use the static method
\verb$domain.disk(a,R)$ which returns this function.
\item The Mesh algorithm can only triangulate simply-connected domains. If your domain is not simply-connected, you need to add cuts to make it simply-connected, else the Mesh algorithm will never terminate.
A cut is typically a segment and can be defined using the static method \verb$domain.segment(a,b)$.
For example, to define a cut in the annulus $R_1<|z|<R_2$, one can do:
\begin{verbatim}
equations=[domain.disk(0,R1),domain.disk(0,R2,complement=True)]
cuts=[domain.segment(0,R2)]
\end{verbatim}
\begin{remark}
There are cases where cuts cannot be segments, for example when using local coordinates (see Section \ref{section:multiDomain}).
So more generally, a cut is defined by a couple $(f,g)$, where $f$ is a real-valued function and $g$ is an assertion (in other words, a boolean-valued function). The cut is defined by $(f(z)=0\text{ and } g(z))$.
The advantage of $g$ being an assertion is that logical operators can be used in its definition.
\end{remark}
\item The parameter \verb$origin$ defines the origin $z_0$ in your domain. This will be the first point of the triangulation, and the integration of your holomorphic data will start there. The origin can be on the boundary of the domain, but should not be on a cut.
\item The parameter \verb$maxNbPoints$ specifies the maximal number of vertices of the triangulation.
This is here to stop Mesh when something gets wrong: triangles are getting to small or the triangulation starts to overlap so will never terminate.
\end{itemize}
\subsection{Defining the holomorphic data}
\label{section:data}
The data is defined by the following lines:
\begin{verbatim}
from lib.cmcR3 import data
mydata=data(potential=eta,Phi0=Phi0)
\end{verbatim}
\begin{itemize}
\item \verb$cmcR3$ is the module which implements the DPW method for CMC surfaces in $\mathbb R^3$.
This must be replaced by the corresponding modules for the other avatars of the DPW method.
\item \verb$eta$ and \verb$Phi0$ are the input data for the DPW method, namely the DPW potential and the initial condition. We will see in Section \ref{section:loops} how loops are implemented in CmcMesh.
For other avatars of the DPW method, the arguments of \verb$data$ may be different but the idea is the same.
\end{itemize}
\subsection{Computing the surface}
\label{section:computingSurface}
The Mesh algorithm is launched by the following line:
\begin{verbatim}
mysurface=mydata.grow(
    mydomain,
    nbdivs=10,
    medDist=0.2,
    medAng=0.1,
    angPref=0.5,
    minCurv=0.01
)
\end{verbatim}
This will at the same time triangulate the domain and compute the corresponding triangulated surface.
The returned value is an instance of the appropriate class, namely \verb$surfaceR3$ for surfaces in euclidean space, \verb$surfaceS3$ for surfaces in the 3-sphere, etc...

The parameter \verb$nbdivs$ specifies the number of subdivisions when integrating numerically the holomorphic data along an edge (using Simpson method, Runge Kutta, etc...).
The parameters \verb$medDist$, \verb$medAng$ and \verb$angPref$ control the length of the edges in the triangulation as follows:
\begin{itemize}
\item \verb$medDist$ is the desired length of edges in the target space,
\item \verb$medAng$ is the desired angle between adjacent normals. This parameter will force the Mesh algorithm to create small triangles where the curvature is large.
\item \verb$angPref$, a real number in the interval $[0,1]$, specifies the preference between the previous two parameters, so \verb$angPref=0$ means that curvature is ignored and \verb$angPref=1$ means that only curvature matters (which can lead to strange results).
\end{itemize}
The formula giving the length of edges issued from a point $p$ in the domain is
$$\mbox{edgeLength}=\frac{1}{\mu}(\mbox{medDist})^{1-\alpha}\times
(\mbox{medAng}/\kappa)^{\alpha}$$
where $\mu$ is the conformal factor of the immersion at $p$ and $\kappa$ is the maximum principal curvature at $p$.
There is a problem in this formula at flat points of minimal surfaces where $\kappa=0$.
The optional parameter \verb$minCurv$ is here for this reason:
$\kappa$ is in fact the maximum of \verb$minCurv$ and the maximum principal curvature, so
$\kappa\geq\mbox{minCurv}$.

The parameters \verb$medDist$, \verb$medAng$ and \verb$angPref$ were already in the historical Mesh algorithm, with the same names. A fine tuning of these parameters yields beautiful triangulations,
especially with the use of multi-domains (see Section \ref{section:multiDomain}).

\subsection{Building the surface}
One can extend the surface by symmetry using the following lines:
\begin{verbatim}
mysurface2=mysurface.replicate(A1,B1,order=2,orient=-1)
fullsurface=mysurface2.replicate(A2,B2,order=3,orient=1)
...
\end{verbatim}
\begin{itemize}
\item For surfaces in $\mathbb R^3$, each $A_i$ is a matrix in the orthogonal group $O(3)$
and $B_i$ is a vector in $\mathbb R^3$, representing the affine isometry $f(x)=A_ix+B_i$.
For other ambiant spaces, the parameters describing the isometry are different but the idea is the same.
\item The parameter \verb$order$ tells how many times the isometry should be applied (usually equal to its order). Namely, if $S$ is your surface, the method \verb$replicate$ returns the union of $f^i(S)$ for
$0\leq i<$\verb$order$.
\item The parameter \verb$orient$ tells wether the isometry preserves (\verb$orient$=1) or reverses
(\verb$orient$=-1) the orientation of the surface. It may happen that an isometry preserves the orientation of the ambiant space but reverses the orientation on the surface. (For example, the order 2 rotation around a line contained in a minimal surface.)
\end{itemize}
Alternately, the surface can be built by providing a list of generators of the isometry group as follows:
\begin{verbatim}
f1=surfaceR3.isometry(A1,B1,orient=-1)
f2=surfaceR3.isometry(A2,B2,orient=1)
...
fullsurface=mysurface.build([f1,f2,...],maxOrder=100)
\end{verbatim}
The \verb$build$ method computes the group generated by the given list of generators and applies it to the surface. It stops when \verb$maxOrder$ elements have been created. It also prints the order of the generated group. See \verb$examples.weierstrass.catenoid$ for a basic example, and 
\verb$examples.cmcS3.lawson_diagonal$ for a more interesting example of isometry group.
\subsection{Curvature lines}
\label{section:curvLines}
On a CMC surface, curvature coordinates can be computed by integrating the square root of the Hopf quadratic differential $Q$:
$$u(z)+i v(z)=\int_{z_0}^z \sqrt{Q}.$$
To compute curvature coordinates, you must set \verb$computeCurvCoords=True$ in the arguments of the \verb$grow$ method. (The default value of this argument is \verb$False$ to shorten computation time.)
You can then extract curvature lines with the method
\begin{verbatim}
mysurface.extractCurvLines(n=5,m=5)
\end{verbatim}
The range of $u$ is divided in $n$ intervals $[a_i,a_{i+1}]$ of the same length and the lines 
$u=a_i$ are computed, for $i=0,\cdots,n$. In the same way, the range of $v$ is divided in $m$ intervals $[b_i,b_{i+1}]$ and the lines $v=b_i$ are computed.
The viewer will display these lines on the surface.

\subsection{Export}
The following lines export the domain and the surface to various formats :
\begin{verbatim}
mydomain.exportMathematica("pictures/domain.txt")
mysurface.exportMathematica("pictures/surface.txt")

mydomain.exportPython("pictures/domain.py")
mysurface.exportPython("pictures/surface.py")

mysurface.exportObj("pictures/surface.obj",scale=5)
\end{verbatim}

\begin{itemize}
\item The Mathematica notebook \verb$view.nb$ allows you to visualize the domain and surface with Mathematica.
\item The python file \verb$view.py$ allows you to visualize the domain and surface using matplotlib.
The matplotlib viewer does not use the GPU so is rather slow. We recommend you only use this viewer when everything else fails.
\item With the Wavefront \verb$OBJ$ format, you can do 3D rendering of your surface, using for example
\href{https://threejs.org/editor/.}{three.js editor}.
If curvature coordinates have been computed (see Section \ref{section:curvLines}), the method will export them as texture coordinates. You can then put a texture in the 3D renderer to display the curvature coordinates. Basic textures are available in \verb$textures$.
\end{itemize}
\subsection{The viewer}
The following lines launch the viewer:
\begin{verbatim}
v=viewer()     
v.addDomain(mydomain,displayBoundary=True)
v.addSurface(mysurface,boundary=mydomain.boundary)
v.addSurface(fullsurface,color='blue',opacity=0.5)
v.show(width=500,height=500)
\end{verbatim}
With these commands, the viewer will display the domain, the fundamental piece and the surface completed by symmetry. You can of course add more commands.
\begin{itemize}
\item The argument \verb$displayBoundary=True$ of \verb$addDomain$ tells the viewer to display the boundary. Each boundary component will be displayed with a different color, depending on the boundary equation it satisfies or the cut it lies on.
\item The argument \verb$boundary=mydomain.boundary$ of \verb$addSurface$ enables the viewer to display the boundary on the surface, with the same colors as in the domain.
\item The argument \verb$color='blue'$ enables you to choose the color from \href{https://examples.vtk.org/site/ColorNamesSeries/ColorNamePatches.html}{Vtk Named Colors}. You can also set \verb$color='normal'$ to color the surface by the normal (as in the historical Mesh),
or \verb$color='curvature'$ to color by the curvature.
\end{itemize}
There are many other optional arguments which are rather intuitive: see the example files.
Here are some functionalities of the viewer:
\begin{itemize}
\item Press 'w' to activate the wireframe mode, so only the edges of the triangulation are visible.
\item Press 's' to return to surface mode.
\item Press 'f' to zoom in where the cursor is.
\item Every click in the domain or the surface will print the closest vertex in the console.
\end{itemize}
\begin{remark}
The viewer can only visualize surfaces in $\mathbb R^3$. If your ambiant space is the 3-sphere or hyperbolic space, you need to project your surface to $\mathbb R^3$ first, using the
\verb$project$ method.
\end{remark}
\section{Implementation of the DPW method}
\subsection{Loops}
\label{section:loops}
In the DPW method, loops are functions defined on the unit circle $\mathbb S^1\subset\mathbb C$.
The variable is usually called $\lambda$.
They are expanded in Fourier series as
$$f(\lambda)=\sum_{k=-\infty}^{\infty} f_k\lambda^k.$$
Theoretically, loops are in some suitable Banach algebra which ensures that the coefficients $f_k$ decay as $k\to\pm\infty$. So loops are implemented by truncating the Fourier series to some large enough order:
$$f(\lambda)=\sum_{k=-n}^{n} f_k\lambda^k.$$
\begin{itemize}
\item The class \verb$loop$ implements complex-valued loops.
The attribute \verb$coeffs$ of a loop is a numpy array containing the $2n+1$ coefficients $f_{-n},\cdots,f_{n}$.
\item The class attribute \verb$loop.deg$ represents the truncation order $n$. It can be changed using the class method
\verb$loop.setDegree$, but this must be done before creating any loop. Once a loop has been created, the truncation order cannot be changed anymore. Default value is 10.
\item The constructor \verb$loop(a,k)$ constructs the loop $a\lambda^k$. For example,
\begin{verbatim}
f=loop(2,-3)+loop(3,0)+loop(1j,4)
\end{verbatim}
constructs the loop $2\lambda^{-3}+3+i\lambda^4$, remembering that \verb$1j$ represents the complex number $i$ in python.
Alternately, you could define
\begin{verbatim}
t=loop(1,1)
f=2*t**-3+3+1j*t**4
\end{verbatim}
Unfortunately, \verb$lambda$ is a reserved keyword of the python language so you cannot define
\verb$lambda=loop(1,1)$. You can, however, use the greek letter $\lambda$ and define
$\lambda$\verb$=loop(1,1)$ but you have to find it on your keyboard or copy it.
\item All arithmetic operations with loops are coded in numpy so are fast. For example, the product of loops uses
\verb$np.convolve(mode='same')$ which implements the truncated convolution product. The class
\verb$loop$ also implements some useful operators for the DPW method:
\begin{itemize}
\item \verb$f.coeff(k)$ returns the coefficient $f_k$.
\item \verb$f.positive()$ returns $f^+=\sum_{k>0} f_k\lambda^k$.
\item \verb$f.negative()$ returns $f^-=\sum_{k<0} f_k\lambda^k$.
\item \verb$f.star()$ returns $f^*=\sum \overline{f_k}\lambda^{-k}$.
\item \verb$f.shift(k)$ multiplies $f$ by $\lambda^k$ by shifting the coefficients.
\item \verb$f.deriv()$ returns the derivative of $f$.
\end{itemize}
\end{itemize}
\subsection{Matrix loops}
\begin{itemize}
\item The class \verb$matrixloop$ implements $\mathcal M(2,\mathbb C)$-valued loops.
It derives from the class \verb$matrix$ of numpy.
For example
\begin{verbatim}
eta=matrixloop([[0,loop(r,-1)+s],[loop(r,1)+s,0]])/z
\end{verbatim}
defines the standard potential for CMC unduloids.
\item The linear ODE $d\Phi=\Phi\eta$ is solved using the Runge Kutta order 4 method, which is implemented as
\verb$solveODE$ in the module \verb$util$.
\item The method \verb$iwasawa()$ returns the Iwasawa factorisation of a loop $\Phi\in\Lambda SL(2,\mathbb C)$
as $\Phi=F\times B$ with $F\in\Lambda SU(2)$ and $B\in\Lambda^+_{\mathbb R}SL(2,\mathbb C)$.
The Iwasawa factorisation algorithm is due to McIntosh-Schmitt. It boils down to solving a $4n\times 4n$
hermitian system using the Cholesky method, where $n=$\verb$loop.deg$ is the truncation order. It is fully implemented in numpy.
\end{itemize}
\section{How the Mesh algorithm work}
\label{section:how}
Each point of the triangulation of the domain is represented by an instance of the class \verb$point$ in the module \verb$mesh$. The main attributes of a point are:
\begin{itemize}
\item \verb$z$, its coordinate in the complex plane,
\item \verb$edgeLength$, which contains the desired length of future edges issued from that point,
\item \verb$Phi$, which contains the integrated holomorphic data from the domain origin to that point.
\end{itemize}
\subsection{Creating edges}
Each avatar (\verb$cmcR3$, \verb$weierstrass$, \verb$triangulate$, etc...) interacts with \verb$mesh$ through the method \verb$nextEdge$ of the class
\verb$domain$.
A call to \verb$nextEdge$ returns a new edge \verb$(p1,p2)$ where \verb$p1$
is an already known point and \verb$p2$ is a new point. It is the responsability of the avatar
to integrate the holomorphic data from \verb$p1$ to \verb$p2$, store the result in \verb$p2.Phi$,
compute the desired edge length from \verb$p2$ as explained in Section \ref{section:computingSurface}
and store it in \verb$p2.edgeLength$.

At each step, the triangulated patch is topologically a disk whose boundary is called the active boundary. When prompted to create a new edge, Mesh looks for the point \verb$p1$ in the active boundary with smallest exterior angle. It divides this angle in a certain number $k$ of angles as close as possible to $\pi/3$ and creates $k$ new triangles. The $k-1$ new points are placed at distance
\verb$p1.edgeLength$ from \verb$p1$.
The $k-1$ new edges are stored in the pending list of edges, and \verb$nextEdge$ returns the first edge in that list.
\subsection{Boundary}
When creating a new edge \verb$(p1,p2)$, Mesh evaluates at \verb$p2$ all equations defining the domain. If an equation $f_i$ returns a negative value, Mesh solves the equation $f_i(z)=0$ by dichotomy
and replaces \verb$p2$ by the solution it founds. Also, the point \verb$p2$ is marked as unactive so no further edge will be issued from it.

Cuts are handled in the same way. If the edge \verb$(p1,p2)$ intersects a cut, Mesh replaces \verb$p2$ by the intersection point and marks it as unactive. Note that as a boundary, each cut has two sides.
\subsection{Corners}
The triangulation is complete when all points on its boundary are marked as unactive.
The final step is to fill in missing corners.
A corner is an intersection point between two boundary components (corresponding to two boundary equations). There is no reason that Mesh has put a point at corners, and this leaves very unpleasant notches.

To fix this, each boundary point has an attribute \verb$bdy$ which contains the index of its boundary component. When two consecutive boundary points \verb$p1$ and \verb$p2$ have different attribute \verb$bdy$, Mesh tries to solve the corresponding two equations using Newton method for a solution \verb$p3$, creates the triangle \verb$(p1,p3,p2)$ and adds the edges \verb$(p1,p3)$ and \verb$(p2,p3)$ to the pending list.

In case Newton method fails to converge, the constructor of the class \verb$domain$ has an optional argument
\verb$corners=[z1,z2,...]$. If defined, Mesh will search that list for a point satisfying the two equations.

When Mesh is done fixing corners, it raises the exception \verb$'done'$ to tell the avatar that the triangulation is complete.

\section{Advanced feature: local coordinates}
\subsection{Motivations}
A new interesting feature of CmcMesh is that it allows you to define your domain using local coordinates and changes of coordinates. In other words, as a manifold.
In all the examples implemented so far, the underlying manifold is in fact a covering of a domain in the plane, so there is no need to introduce the manifold structure. But local coordinates are still interesting for the following reasons:
\begin{itemize}
\item Without local coordinates, Mesh cannot triangulate the Riemann sphere $\mathbb C\cup\{\infty\}$: however fast the triangles grow, the algorithm will never terminate. Using the local coordinate $w=\frac{1}{z}$ in the chart $|z|>1$, it succeeds in triangulating the Riemann sphere, with no change in the algorithm. See \verb$examples/triangulate/riemann_sphere$.
\item Local coordinates seem to be the right answer to the so-called shearing problem.
Consider the module \verb$examples.weierstrass.scherk$ which implements the classical Scherk minimal surface. On the left image of Figure \ref{fig:scherk}, you can see that the triangulation becomes very irregular when approaching the end. The reason is that the immersion is asymptotic to
$(\log|z-p|,0,\arg(z-p))$ as $z$ tends to the puncture $p$. And this map, even though it is conformal, distorts triangles very much near $p$.
It is clear that $w=\log(z-p)$ is a more appropriate coordinate in a neighborhood of $p$.
Using $w=\log(z-p)$ as a local coordinate, we obtain a much nicer triangulation: see the right image of Figure \ref{fig:scherk} and the module \verb$examples.weierstrass.scherk2$.
\begin{figure}
\includegraphics[height=2cm]{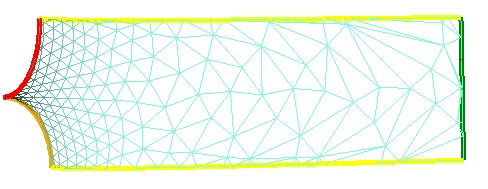}
\hspace{1cm}
\includegraphics[height=2cm]{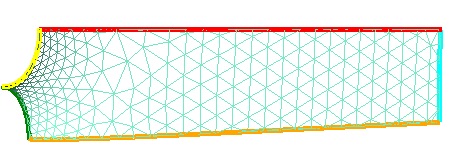}
\label{fig:scherk}
\caption{Left: triangulation of the Scherk surface without using local coordinates.
Right: triangulation of the Scherk surface using the local coordinate $w=\log(z-p)$
in a neighborhood of the puncture $p=e^{i\pi/4}$.}
\end{figure}
\item In the DPW method, examples often have apparent singularities, where the potential has a pole but can be gauged to a smooth potential using a local regularizing gauge.
It is often the case that at apparent singularities, a change of coordinate is also necessary.
See the example \verb$advanced.cmcS3.lawson_genus2$ for an illustration of the use of local coordinates and regularizing gauges.
\end{itemize}
\subsection{Definition of a multi-domain}
\label{section:multiDomain}
The class \verb$multiDomain$ allows you to define a domain as a manifold.
This is in fact the fundamental class of the module \verb$mesh$: a domain is a multi-domain with only one chart and the class \verb$domain$ derives from \verb$multiDomain$.
A multi-domain is defined with the following lines:
\begin{verbatim}
from lib.mesh import multiDomain
mydomain=multiDomain(
   equations=[f1,f2,...],
   cuts=[cut1,cuts2,...],
   origin=z0,
   originChart=i0,
   maxNbPoints=1000,
   getChart=[c1,c2,...],
   changeCoord=mychangecoord)
\end{verbatim}
\begin{itemize}
\item The argument \verb$getChart$ is a list of $n$ integer-valued functions, and
\verb$changeCoord$ is an $n\times n$ matrix of complex-valued functions, where
$n$ is the number of charts. They work as follows:

Let us say our manifold is covered by $n$ charts $U_0,\cdots,U_{n-1}$, with local complex coordinates $\varphi_i:U_i\to\mathbb C$.
Each point has a prefered chart, given by the function \verb$getChart$, and changes of coordinates
are implemented by the function \verb$changeCoord$.
More precisely, if $p\in U_i$ has coordinate $\varphi_i(p)=z$, the lines
\begin{verbatim}
j=getChart[i](z)
w=changeCoord[j][i](z)
\end{verbatim}
return the integer $j$ such that $U_j$ is the prefered chart of $p$, 
and the complex coordinate $w=\varphi_j(p)$ of $p$ in this chart.
Of course, if $U_i$ is the prefered chart of $p$, then $j=i$ and $w=z$.

See the module \verb$examples.triangulate.riemann_sphere$ for the definition of the Riemann sphere using its standard atlas with two charts.
\item The arguments \verb$origin$ and \verb$originChart$ allow you to specify the origin in your multi-domain.
The default value of \verb$originChart$ is zero.
\item Each member of the list \verb$equations$ is a list of $n$ equations, representing the equation of the boundary component in each chart. In other words, a point $p\in U_i$ with coordinate $z=\varphi_i(p)$
is in the domain if for all $k$, \verb$equations[k][i](z)$$>0$.
If the corresponding boundary component does not intersect the chart $U_i$, you can set
\verb$equations[k][i]=Null$.

\item In the same way, each member of the list \verb$cuts$ is a list of $n$ couples defining the cut in each chart. If the cut does not intersect the chart $U_i$, you can set
\verb$cuts[k][i]=Null$.
\end{itemize}
See \verb$examples.triangulate.riemann_sphere_with_holes$ for the definition of a multi-domain using two charts, several boundary equations and several cuts. In this example, the change of coordinate is $z\mapsto 1/z$, so a cut which is a segment in one chart may be a circle in the other chart.
\subsection{Definition of the holomorphic data on a multi-domain}
The holomorphic data on a multi-domain is defined using the following lines:
\begin{verbatim}
from lib.cmcR3 import multiData
mydata=multiData(
   potential=[eta0,eta1,...],
   Phi0=Phi0,
   gauge=[G0,G1,...])
\end{verbatim}
\begin{itemize}
\item The module \verb$lib.cmcR3$ implements the DPW method for CMC surfaces in $\mathbb R^3$.
Each avatar of the DPW method has its own class \verb$multiData$.
As already said, in the DPW method, local coordinates work in conjonction with regularizing gauges to deal with apparent singularities. Therefore, we assume that each chart $U_i$ comes with a regularizing gauge $G_i$.
\item The argument \verb$gauge$ contains the list of the gauges $G_i$, expressed in the local coordinate of the chart $U_i$. If no gauge is needed, you may set
$G_i=$\verb$lambda z:matrixloop([[1,0],[0,1]])$.
\item The argument \verb$eta$ contains the list of the gauged potentials, expressed in local coordinates.
In other words, $\eta_i = (\varphi_i^*\eta)\cdot G_i$.
\end{itemize}
See the module \verb$examples.cmcR3.sphere2$ for a basic example of the use of local coordinates and regularizing gauges, and
\verb$advanced.cmcS3.lawson_genus2$ for a more elaborate example.
\subsection{Visualization of multi-domains}
The method \verb$addDomain$ has an optional argument \verb$chart$, with default value $0$.
\begin{itemize}
\item With \verb$chart=i$, it displays the chart number $i$.
\item With \verb$chart=[i1,i2,...]$, it displays charts $i_1$, $i_2$, $\cdots$, converting all coordinates to chart $i_1$ using the changes of coordinates.
\end{itemize}
The method must be called several times to display all charts separately.

\bibliography{bib}{}
\bibliographystyle{plain}

\end{document}